\documentclass[12pt,letterpaper]{amsart} 

\usepackage{times}
\usepackage{amsmath}
\usepackage[dvips]{graphics,epsfig}

\newcommand{\Area}{\operatorname{Area}}
\newcommand{\Wr}{\operatorname{Wr}}
\newcommand{\Lk}{\operatorname{Lk}}
\newcommand{\Tw}{\operatorname{Tw}}

\newcommand{\R}{\mathbf{R}}
\newcommand{\cross}{\times}

\newtheorem{theorem}{Theorem}
\newtheorem{definition}[theorem]{Definition}
\newtheorem{prop}[theorem]{Proposition}
\newtheorem{corollary}[theorem]{Corollary}
\newtheorem{problem}[theorem]{Problem}
\newtheorem{conjecture}[theorem]{Conjecture}
\newtheorem{proposition}[theorem]{Proposition}
\newtheorem{lemma}[theorem]{Lemma}

\newcommand{\Calugareanu}{C\u{a}lug\u{a}reanu}

\title[On comparing the writhe of a smooth curve to the writhe of an inscribed polygon]
{On comparing the writhe of a smooth curve \\ to the writhe of an inscribed polygon} 
\author{Jason Cantarella}

\begin{document}
\begin{abstract}
We find bounds on the difference between the writhing numbers of a
smooth curve and a polygonal curve inscribed within. The proof is
based on an extension of Fuller's difference of writhe formula to the
case of polygonal curves. The results establish error bounds useful
in the numerical computation of writhe.
\end{abstract}

\maketitle

\section{Introduction \label{sec:intro}}
The writhing number measures the wrapping and coiling of space
curves. Writhe has proved useful in molecular biology, where it is
used to study the geometry of tangled strands of DNA \cite{sumners};
often with the famous \Calugareanu -White formula for a curve
$C$ in space with a normal field~$V$ \cite{Cal1,Cal2,white,pohl}:
\begin{equation*}
\Lk(C,C+\epsilon V) = \Tw(C,V) + \Wr(C). 
\end{equation*}
In these applications, and in numerical simulations performed by
biologists and mathematicians, it is often required to compute
writhing numbers using numerical methods. 

Several authors have presented algorithms for computing the
exact writhing number of an $n$-edge polygonal curve in a finite
number of steps \cite{cimasoni,banchoff,edelsbrunner,sumners}. The
fastest of these algorithms runs in time $O(n^{1.6})$, while earlier
methods use time $O(n^2)$.

Careful implementations of such algorithms provide acceptable accuracy
in computing writhe for polygonal curves. But reliably computing the
writhe of smooth curves requires another step: we must be able to
bound the error introduced in approximating a smooth curve by an
inscribed polygonal curve. The purpose of this paper is to prove:

\begin{theorem} \label{thm:main}
Suppose $C(t)$ is a simple, closed curve of class $\mathcal{C}^4$.  We
assume $C(t)$ is parametrized so that $|C'(t)| \geq 1$, and that
we have upper bounds $B_1, \dots, B_4$ on $|C'(t)|, \dots
|C^{(4)}(t)|$.  Let $C_n(t)$ be any $n$-edge polygonal curve inscribed
in $C$ with maximum edge length $x$ and $1/x > 5B_2$.

If the ribbon formed by joining $C_n(t)$ to $C(t)$ for every $t$ is
embedded,
\begin{equation}
| \Wr(C) - \Wr(C_n) | < \alpha \, n x^3 + n O(x^4).
\end{equation} 
where $\alpha$ is a numerical constant less than $B_2 (5 B_2^2 + B_3)$.
\end{theorem}
That is, if the lengths of the edges of $C_n$ are approximately constant,
the error is bounded by a multiple of $1/n^2$. 

The proof is based on Fuller's $\Delta\!\Wr$ formula, which gives the
difference in writhing number between two curves as the spherical area
of the ribbon bounded by the curves on $S^2$ swept out by their unit
tangent vectors \cite{fuller}. (Following Bruce Solomon
\cite{solomon}, we will refer to such curves as {\em tantrices},
though they are classically referred to as {\em tangent
indicatrices}.)

We begin by defining the writhing number in
Section~\ref{sec:definitions}. Sections~\ref{sec:fuller}
and~\ref{sec:area} then introduce the original form of Fuller's
$\Delta\!\Wr$ formula.  In Sections~\ref{sec:extension_I} and
\ref{sec:extension_II} we extend Fuller's formula to the case where
one curve is polygonal and the other is of class $\mathcal{C}^2$ using a
natural geometric idea: the tantrix of a polygonal curve should be
defined to be the chain of geodesic segments on $S^2$ joining the
(isolated) tangent vectors of the curve (this was pointed out by Chern
in \cite{Chern}). In the process, we discover a surprising fact: the
writhe of a polygonal curve is {\em equal} to the writhe of any smooth
curve obtained by carefully rounding off its corners!

Section~\ref{sec:estimate} contains the remainder of our work:
estimating the terms in our improved version of the $\Delta\!\Wr$
formula to obtain Theorem~\ref{thm:disc-error}. We test our error 
bounds in Section~\ref{sec:examples} by computing the writhe of 
a collection of polygonal curves inscribed in a smooth curve of 
known writhe. 

The last section contains a discussion of some open problems inspired
by the present work. We state the most important of them now: Like
most of the theory of writhing numbers, the proof of our main theorem
depends essentially on the fact that $C$ is closed. Can these
methods be extended to open curves?

\section{Definitions\label{sec:definitions}}

The writhing number of a space curve is defined by:
\begin{definition} \label{def:writhe}
The {\em writhe} of a piecewise differentiable curve $C(s)$ is
given by:
\begin{equation} \label{eqn:writhe}
\Wr(C) = \frac{1}{4 \pi} \int_{C \cross C} \frac{ C'(s) \cross
C'(t) \cdot (C(s) - C(t)) } { |C(s) - C(t)|^3 } \,\mathrm{d}s\,\mathrm{d}t,
\end{equation}
\end{definition}

Definition~\ref{def:writhe} is inspired by the Gauss formula for the
linking number of two space curves, $A(s)$ and $B(s)$ (see Epple
\cite{epple} for a fascinating discussion of the history of this formula):
\begin{equation}
\Lk(A,B) = \frac{1}{4 \pi} \int_{A \cross B} \frac{ A'(s) \cross
B'(t) \cdot (A(s) - B(t)) } { |A(s) - B(t)|^3 } \,\mathrm{d}s\,\mathrm{d}t.
\end{equation}

When the two curves $A$ and $B$ become a single curve, their linking
number becomes the writhing number. This introduces a potential
singularity on the diagonal of $C \cross C$, but a careful calculation
shows that the integral still converges.  In fact, the integrand of
Equation~\ref{eqn:writhe} approaches $0$ on the diagonal of $C \cross
C$, even when the curve $C$ has a corner. 

From now on, we'll assume that $C$ is simple. With this assumption,
another way to look at the integral of Definition~\ref{def:writhe} is
to observe that the integrand is the pullback of the area form on
$S^2$ under the Gauss map $C \cross C \rightarrow S^2$ defined by
\begin{equation}
(C(s),C(t)) \mapsto \frac{C(s) - C(t)}{|C(s) - C(t)|}.
\end{equation}
From this perspective, we can see that the (signed) multiplicity of
the Gauss map at any point~$p$ on~$S^2$ is just the number of
self-crossings of the projection of~$C$ in direction~$p$.  

\section{Fuller's $\Delta\!\Wr$ formula \label{sec:fuller}}

Suppose we have a differentiable curve $C(t)$, with unit tangent
vector $T(t)$.  As we mentioned in Section~\ref{sec:intro}, the curve
$T(t)$ on the unit sphere is known as the tantrix of $C$. This curve
divides the unit sphere into a number of cells. Within each cell, the
signed crossing number of the projection of~$C$ is constant: changing
projection directions within the cell amounts to altering the
projection of the knot by a regular isotopy consisting of Reidemeister
moves of type II and III (pictured below). Neither of these moves
changes the signed crossing number of the knot.

\medskip

\begin{center}
\includegraphics{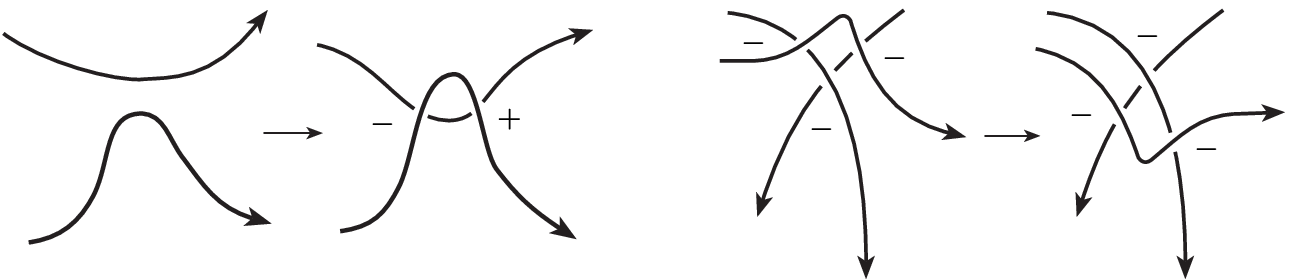}
\end{center}
\smallskip
\begin{center}
\parbox{5in}{\noindent{\bf Figure 1.} Changing the projection direction within a cell can only
alter the diagram by one of these two moves. Neither changes the signed
crossing number of the diagram, as we can see by counting the $+$ and $-$
markers at the crossings of $C$.}
\end{center}

\medskip

This observation motivates the idea that the writhe of a closed space
curve is related to the fraction of the sphere's area enclosed by its
tantrix. In 1978, Brock Fuller stated the following:

\begin{theorem} {\rm (Fuller's Spherical Area Formula)} \label{thm:spherical_area}
For any closed space curve $C(s)$ of class $\mathcal{C}^3$, let $A$ be
the spherical area enclosed by the tantrix of $C$. Then
\begin{equation}
1 + \Wr(C) = \frac{A}{2\pi} \!\! \mod 2.
\end{equation}
\end{theorem}

Fuller used this formula to conclude that the difference in writhe
between two curves $X_0$ and $X_1$ whose tantrices $T_0$ and $T_1$ are
sufficiently close is given by a certain formula, which represents the
spherical area of the ribbon between $T_0$ and $T_1$.

To be more specific, suppose that $X_0$ and $X_1$ are simple closed
space curves of class $\mathcal{C}^2$, with regular parametrization
(that is, parametrized so that $X_0'$ and $X_1'$ never vanish), and
unit tangent vectors $T_0$ and $T_1$. Let $F \!:\! S^1 \cross [0,1]
\rightarrow \R^3$ be a continuous deformation of $X_0$ into $X_1$,
where $F(t,\lambda) = X_\lambda(t)$ and the $X_\lambda$ are simple
curves of class $\mathcal{C}^1$, with unit tangent vectors
$T_\lambda(t)$ continuous in $(t,\lambda)$.

\begin{theorem}{\rm (Fuller's $\Delta\!\Wr$ Formula)} 
\label{thm:fuller_delta_writhe}
If $T_1(t)$ and $T_\lambda(t)$ are not antipodal for all
$(t,\lambda)$, then 
\begin{equation} \label{eqn:fuller_diff}
\Wr(X_1) - \Wr(X_0) = \frac{1}{2\pi} \int_C 
\frac{T_0(t) \cross T_1(t)}{1 + T_0(t) \cdot T_1(t)} \cdot
\left[ T'_0(t) + T'_1(t) \right] \, \mathrm{d}t.
\end{equation}
\end{theorem}

We observe that this formula does not require an arc-length
parametrization of $X_0$ and $X_1$.

\section{Justifying Fuller's interpretion of the $\Delta\!\Wr$ formula \label{sec:area}}

While Fuller stated both these theorems in 1978, he did not provide
complete proofs for either. The first rigorous proofs of
Theorems~\ref{thm:spherical_area} and~\ref{thm:fuller_delta_writhe} were given
by Aldinger, Tabor, and Klapper \cite{ATK} in~1995. While these
authors proved both theorems as stated, they did not show that the
formula in Theorem~\ref{thm:fuller_delta_writhe} represents the spherical area
of the ribbon between $T_0$ and $T_1$ (in \cite{ATK}, the right-hand
side of Equation~\ref{eqn:fuller_diff} describes the difference
between the twist of two frames on $X_0$ and $X_1$.) 

In the spirit of their paper, we now justify Fuller's original
intuition about Equation~\ref{eqn:fuller_diff}.

\begin{prop} \label{prop:ribbon}
Given two curves $T_0(t), T_1(t)\!:\! [0,1] \rightarrow S^2$ where
$T_0(t)$ and $T_1(t)$ are never antipodal, the area of the spherical
region $R$ bounded by $T_0$, $T_1$ and the great circle arcs joining
their endpoints is given by
\begin{equation}
\Area(R) = \int \frac{T_0(t) \cross T_1(t)}{1 + T_0(t) \cdot T_1(t)} \cdot \left(T'_0 + T'_1 \right) \, \mathrm{d}t.
\end{equation}
\end{prop}

\begin{proof}
We let 
\begin{equation*}
u(\theta,t) = \cos \theta \,\,  T_0(t) + \sin \theta \,\, T_1(t),
\end{equation*}
and parametrize the region $R$ by
\begin{equation*}
v(\theta,t) = \frac{u(\theta,t)}{|u(\theta,t)|}
\end{equation*}
where $\theta$ ranges from $0$ to $\pi/2$. Plugging this
parametrization into the area form on $S^2$, and using the properties
of the triple product, we find
\begin{equation*}
\mathrm{d}\Area = \frac{1}{|u|^3} 
	\left( 
	\frac{\partial u}{\partial \theta} \cross 
    	\frac{\partial u}{\partial s} \cdot u
	\right) \mathrm{d}\theta \wedge \mathrm{d}t.
\end{equation*}
Using the definition of $u(\theta,t)$, this simplifies to 
\begin{equation*}
\mathrm{d}\Area = T_0 \cross T_1 \cdot 
	\left( 
	\frac{\cos \theta}{(1 + 2 \cos \theta \sin \theta \,\, T_0 \cdot T_1)^{\frac{3}{2}}} T'_0 +  
	\frac{\sin \theta}{(1 + 2 \cos \theta \sin \theta \,\, T_0 \cdot T_1)^{\frac{3}{2}}} T'_1
	\right) \mathrm{d}\theta \wedge \mathrm{d}t.
\end{equation*}
Using the formula $\sin 2\theta = 2 \cos\theta \sin\theta$, and the
fact that the definite integrals of each of the trigonometric
expressions above from $0$ to $\pi/2$ are equal, we have
\begin{equation*}
\Area(R) = \int_0^1 T_0 \cross T_1 \cdot
	\left[ 
	\int_0^{\pi/2} \frac{\cos \theta}{(1 + \sin 2 \theta \,\,T_0 \cdot T_1)^{3/2}} \, \mathrm{d}\theta 
	\right] (T'_0 + T'_1) \, \mathrm{d}t.
\end{equation*}
This can be solved by the general integration formula
\begin{equation}
\int \frac{\cos \theta}{(1 + a \sin 2\theta)^{3/2}} \, \mathrm{d}\theta =
\frac{-a \cos \theta - \sin \theta}{(a^2 - 1) \sqrt{1 + a\sin 2\theta}},
\end{equation}
which yields the formula in the statement of the Proposition.
\end{proof}

\section{Extending Fuller's formula to polygonal curves: I \label{sec:extension_I}}

To measure the difference in writhe between a smooth curve and a
polygonal curve inscribed in the smooth curve, we must extend
Theorem~\ref{thm:fuller_delta_writhe} to polygonal curves.  To do so,
we intend to approximate each polygonal curve with a family of smooth
curves so that the writhe of the smooth curves converges to the writhe
of the polygonal curve. 

Examining Definition~\ref{def:writhe}, it might seem that this result
follows from general principles. For instance, one might conjecture
that $\Wr$ was continuous in the $C^1$ norm on curves, and hope to
obtain an approximating family using standard
techniques. Unfortunately, the situation is not so simple; as the
example in Figure~2 shows, writhe is not continuous in any $C^k$ norm
on curves. Thus, our proof depends explicitly on the hypothesis that
the limit curve is polygonal; it cannot be easily extended to the case
where the limit curve is merely piecewise $C^2$.

\begin{center}
\includegraphics{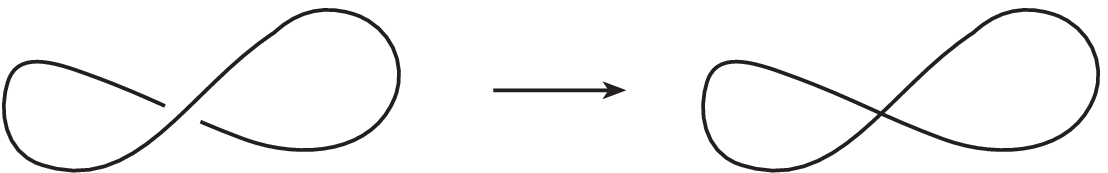}
\end{center}
\smallskip
\begin{center}
\parbox{5in}{\noindent{\bf Figure 2.} The family of almost-planar curves on the
left converge in any $C^k$ norm to the planar figure eight curve on
the right. However, the writhe of the curves on the left approaches
one, while the writhe of the planar figure eight is zero. This shows
that writhe is not continuous in any $C^k$ norm on curves.}
\end{center}
\medskip

To prepare for the proof, we establish some notation for polygonal
curves. Let $C(t)$ be a polygonal curve with corners at cyclically
ordered parameter values $t_0 < t_1 < \dots < t_n = t_0$. We
let $T(t)$ denote the unit tangent to $C$, and set up the convention
that $T(t_i)$ will be the tangent vector {\em leaving} $C(t_i)$.

We now construct a family of smooth curves approximating our polygonal curve. 

\begin{prop} \label{prop:approx}
Given an embedded polygonal curve $C$ with corners at $t_0, \dots, t_{n-1}, t_{n} = t_0$,
there exists a family of smooth curves $C_i$ converging pointwise to
$C$ with
\begin{enumerate}
\item $C_i = C$ outside a neighborhood of each corner point $C(t_j)$
of radius $1/i$.
\item Near each corner, the tangent vectors of $C_i$
      interpolate between $T(t_{j-1})$ and $T(t_j)$.
\item $\Wr(C_i) \rightarrow \Wr(C)$.
\end{enumerate}
\end{prop}

\begin{proof}
It is easy to construct a family of $C_i \rightarrow C$ obeying
conditions (1) and (2) by rounding off each corner of $C$. We claim
that this can be done in such a way that the writhe integrand has a
uniform upper bound on all the $C_i$. Since condition (1) implies that
the $C_i \rightarrow C$ pointwise in the~$\mathcal{C}^1$ norm, the
bounded convergence theorem~\cite[p.81]{royden} will then yield
condition~(3).

Since any pair of adjacent edges is planar, we can choose the $C_i$ so
that the region of each $C_i$ approximating a pair of adjacent edges
is also planar. This means that for some universal $\epsilon$, the
writhe integrand of each $C_i$ vanishes in an $\epsilon$-neighborhood
of the diagonal of $C_i \cross C_i$.

Since $C$ has no self-intersections and the angle at each corner of
$C_i$ is positive, the distance between any pair of non-adjacent edges
of $C$ is bounded below by some constant. Since the $C_i$ converge to
$C$ pointwise, we may assume the same for the portions of the $C_i$
approximating any pair of disjoint edges. Throwing away finitely many
of the $C_i$ if necessary, this means that for any $\delta > 0$, there
exists a universal lower bound (depending on $\delta$) on the distance
between any pair of points in $C_i \cross C_i$ {\em outside} an
$\delta$-neighborhood of the diagonal.

But for any pair of points on $C_i$, the writhe integrand is bounded
above by the inverse square of the distance between them. Thus, our
lower bound on self-distances yields a universal upper bound on the
writhe integrand for $C$ and all the $C_i$ outside a
$\delta$-neighborhood of the diagonal. Choosing $\delta < \epsilon$,
this completes the proof of the proposition.
\end{proof}

\section{Extending Fuller's formula to polygonal curves: II \label{sec:extension_II}}

We now state our extension of Fuller's theorem.  Our formula will
apply to the following situation (c.f. Section~\ref{sec:fuller}):
Suppose that $X_0$ and $X_1$ are simple closed space curves, with
$X_0$ of class~$\mathcal{C}^2$ and $X_1$ polygonal, with regular
parametrization (that is, parametrized so that $X_0'$ and $X_1'$ never
vanish where they are defined), and unit tangent vectors $T_0$
and~$T_1$.

Let $F\!:\! S^1 \cross [0,1] \rightarrow \R^3$ be a $C^0$ deformation of
$X_0$ into $X_1$, where $F(t,\lambda) = X_\lambda(t)$, so that the
$X_\lambda$ are simple curves of class $C^1$ for $\lambda \in [0,1)$,
with unit tangent vectors $T_\lambda(t)$ continuous in~$(t,\lambda)$.
As above, we take the corners of $X_1$ to be at parameter
values $t_0, t_1 , \dots , t_n = t_0$.  We let $T_1$ denote
the unit tangent vector to $X_1$, and let $T_1(t_i)$ be the tangent
vector {\em leaving} $X_1(t_i)$.

\begin{theorem} \label{thm:extended_fuller} 
If each corner angle of $X_1$ is strictly greater than $\pi/2$, and
each $T_1(t)$ and $T_\lambda(t)$ are at an angle less than $\pi/2$,
then
\begin{equation*} \label{eqn:diff}
\Wr(X_1) - \Wr(X_0) = 
\frac{1}{2\pi} \sum_{i=1}^{n}
  \Area R(T_0(t_i),T_0(t_{i+1}),T_1(t_i)) 
+ \Area \triangle T_0(t_i) T_1(t_{i-1}) T_1(t_i),
\end{equation*}
where $R(T_0(t_i),T_0(t_{i+1}),T_1(t_i))$ is the spherical region
bounded by geodesics from $T_1(t_i)$ to $T_0(t_i)$ and
$T_0(t_{i+1})$ and the portion of $T_0$ between $t_i$ and $t_{i+1}$,
$\triangle T_0(t_i) T_1(t_{i-1}) T_1(t_i)$ is the spherical
triangle with these three vertices, and $\Area$ represents oriented
area on $S^2$.
\end{theorem}

\begin{proof}
Construct a sequence of smooth curves $C_j \rightarrow X_1$ using
Proposition~\ref{prop:approx}. For large enough~$j$, each of these
curves can be homotoped to $X_1$ through a family of simple
$\mathcal{C}^1$ curves with a continuous family of tangent vectors, as
in the setup for the statement of this theorem above.

Joining these homotopies to the homotopy from $X_1$ to $X_0$ assumed
by our hypotheses generates a family of (non-smooth) homotopies from
the $X_0$ to each of the $C_j$. We wish to smooth each of these to
obtain homotopies from $X_0$ to $C_j$ which obey the conditions of
Fuller's $\Delta\!\Wr$ formula
(Theorem~\ref{thm:fuller_delta_writhe}).

We first prove that the tangent vectors of each of the intermediate
curves in each homotopy from~$X_0$ to $C_j$ are never antipodal to the
corresponding tangent vectors $T_j$ of $C_j$. By hypothesis, for
each~$t$ and~$\lambda$, $\angle T_\lambda(t), T_1(t) < \pi/2$.  On the
other hand, since the difference between the tangent vectors to $X_1$
at any corner is less than $\pi/2$, for large enough $j$, $\angle
T_1(t), T_j(t) < \pi/2$. Putting these equations together, we see that
$\angle T_\lambda(t), T_j(t) < \pi,$ and so these vectors are never
antipodal. 

It is easy to smooth the combined homotopy from $X_0$ to $C_j$ so that
each of the intermediate curves is of class $\mathcal{C}^1$ while
preserving this condition. Since the smoothed homotopy satisfies the
hypotheses of Fuller's $\Delta\!\Wr$ formula
(Theorem~\ref{thm:fuller_delta_writhe}), Proposition~\ref{prop:ribbon}
tells us that the difference between $\Wr(X_0)$ and $\Wr(C_j)$ is the
spherical area of the ribbon joining~$T_0$ and~$T_j$.

For each $i$, the contribution to the spherical area from the straight
part of $C_j$ between $t_i$ and $t_{i+1}$ comes from the ribbon
between $T_1(t_i)$ and the portion of $T_0$ with $t \in (t_i +
1/j,t_{i+1} - 1/j)$. As $j \to \infty$, this area converges to the
area of the ribbon between the portion of $T_0$ with $t \in
(t_i,t_{i+1})$ and $T_1(t_i)$. This is the first term in our sum
above.

At each vertex $t_i$ of $X_1$, the contribution to our spherical area
from the curved part of $C_j$ comes from the ribbon between the great
circle arc connecting $T_1(t_{i-1})$ and $T_1(t_i)$ and a portion of
$T_0$ of parameter length $2/j$. As $j \to \infty$, the area of this
ribbon converges to the area of the spherical triangle with vertices
$T_0(t_i)$, $T_1(t_i)$, $T_1(t_{i-1})$. This is the second term in our
sum above. Figure~3 shows both these terms on the unit sphere.

\begin{center}
\includegraphics{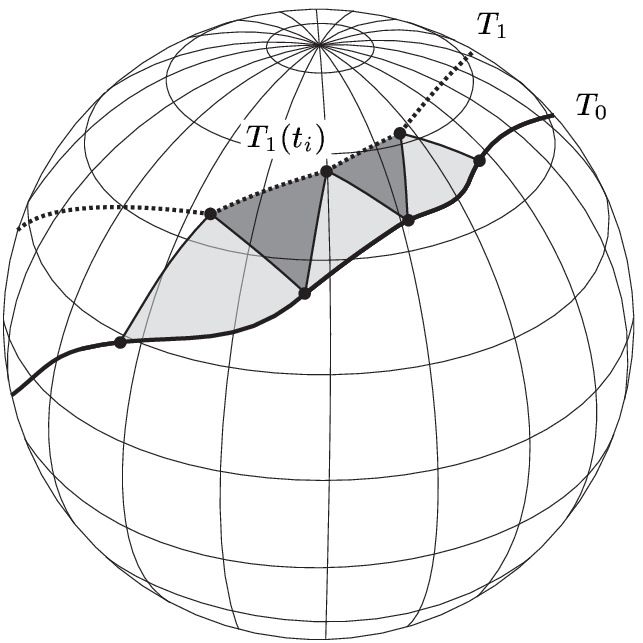}
\end{center}
\smallskip
\begin{center}
\parbox{5in}{
\noindent{\bf Figure 3.} This figure shows the two types of regions in the sum
in the statement of Theorem~\ref{thm:extended_fuller}. The top
(dotted) curve shows the great circle arcs joining the tangent vectors
$T(t_i)$ of the polygonal curve $X_1$. The bottom curve shows the
continuous curve of unit tangents $T_0$ to the smooth curve $X_0$.
The light gray regions show the first terms in the sum, while the dark
gray spherical triangles show the second terms.}
\end{center}
\medskip

We have shown that the right-hand side of the statement of the Theorem
is equal to the limit $\lim_{j \to \infty} (\Wr(C_j) - \Wr(X_0))$.
However, by Proposition~\ref{prop:approx}, $\lim_{j \to \infty}
\Wr(C_j) = \Wr(X_1)$.  Thus
\begin{equation}
\lim_{j \to \infty} \Wr(C_j) - \Wr(X_0) = \Wr(X_1) - \Wr(X_0),
\end{equation}
which is the left-hand side in the statement of the Theorem. This
completes the proof.
\end{proof}

We now make a surprising observation: Since the tantrices of the $C_j$
differ as curves on $S^2$ only in parametrization, the area between
each of these curves and the tantrix of $X_0$ is constant. Thus, by
Fuller's formula, each $C_j$ has the same writhe! And since (by
Proposition~\ref{prop:approx}) these writhing numbers converge to the
writhe of $X_1$, each $\Wr(C_j)$ is equal to $\Wr(X_1)$ as well! So we
have the following corollary:

\begin{corollary}
If $C_n$ is a polygonal curve, and $C$ is a smooth curve obtained by
rounding off the corners of $C_n$ under the conditions of
Proposition~\ref{prop:approx}, then
\begin{equation}
\Wr(C_n) = \Wr(C).
\end{equation}
\end{corollary}

\section{Bounding the $\Delta\!\Wr$ Formula \label{sec:estimate}}

We now prove our main theorem by finding asymptotic bounds for
Fuller's~$\Delta\!\Wr$ formula. Our theorem deals with the
following situation: Assume that $C(t)$ is a simple closed curve of
class~$\mathcal{C}^4$, parametrized so that $|C'(t)| \geq 1$.  (Given
any initial parametrization, this can be accomplished by rescaling.)
Further, assume we have upper bounds $B_1$, \dots, $B_4$ on the norms
of the first four derivatives of $C$.

Let $C_n(t)$ be any $n$-edge polygonal curve inscribed in~$C$. 
We assume that the maximum edge length of $C$ is bounded by~$x$.

\begin{theorem} \label{thm:disc-error}
If the ribbon formed by joining $C_n(t)$ to $C(t)$ for every $t$ is
embedded, and $1/x > 5B_2$,
\begin{equation}
| \Wr(C) - \Wr(C_n) | < \alpha \, n x^3 + n O(x^4).
\end{equation} 
where $\alpha$ is a numerical constant less than $B_2 (5 B_2^2 + B_3)$.
\end{theorem}

We make a few comments on this theorem before diving into the
proof. First, we observe that if the lengths of the edges of $C_n$ are
all of the same order of magnitude, the difference between the writhe
of $C$ and the writhe of $C_n$ is of order $1/n^2$.

Next, we observe that the form of our theorem was chosen to be of
maximal use in applications. In particular, we did not require that
$C$ be parametrized by arclength and state our bounds in terms of
curvature and torsion because in practice it is very difficult
to obtain an arc-length parametrization of a given curve, while it is
comparatively easy to obtain values for the derivative bounds given
above.

Last, we discuss the role of the additional hypotheses in the
statement above; that the ribbon between $C$ and $C_n$ be embedded and
that $1/x$ be greater than $5B_2$.  Both are intended to exert enough
control over the approximation to guarantee the existence of a
homotopy from $C$ to $C_n$ obeying the requirements of
Theorem~\ref{thm:extended_fuller}.

We can guarantee that $C_n$ satisfies the first hypothesis by proving
that $C_n$ lies in an embedded tubular neighborhood of $C$. Since $C$
is of class $\mathcal{C}^4$, and has no self-intersections, such a
neighborhood is guaranteed to exist: for a discussion of how to
compute the radius of this tube (which is known as the {\em thickness}
of $C$), see the literature on {\em ropelength} of knots (\emph{e.g.}
\cite{GM,CKS,LSDR}).

\begin{proof} 
We begin by reparametrizing our curve by arclength. This forces us to
recompute our bounds for the derivatives of $C(t)$ (a standard
computation), arriving at
\begin{equation}
|C'(s)| = 1, \qquad |C''(s)| < K := 2 B_2 , \qquad |C'''(s)| < T := 2 B_3 + 10 B_2^2,
\end{equation}
while $C^{(4)}(s)$ is again bounded above. To remind ourselves of the
connection between these bounds and the curvature and torsion of our
curve, we will refer to the bound for the second derivative as $K$,
and the bound for the third derivative as $T$.  Further, we note that
the curvature $\kappa(s)$ of our curve is bounded above by $K$, and
that our hypotheses imply that $1/x > (5/2) K$.

We also establish the convention that the corners of $C_n$ are at parameter
values cyclically ordered as $s_0, \dots, s_{n-1}, s_n = s_0$. 

By smoothing the linear interpolation between $C$ and $C_n$, we can
construct a homotopy between $C$ and $C_n$ according to the conditions
of Theorem~\ref{thm:extended_fuller} as long as:
\begin{enumerate}
\item the ribbon joining $C$ to $C_n$ is embedded,
\item the angle at each corner of $C_n$ is at least $\pi/2$,
\item the angle between $T(s)$ and $T_n(s)$ is at most $\pi/2$ for any $s$.
\end{enumerate}

Borrowing from Lemma~\ref{lem:angle} (below), we see that our
assumption that $1/x > (5/2) K$ is enough to bound the angle in (3) by
$0.20402 < \pi/4$.  At any corner $s_i$, the same Lemma implies that
the corner angle is the supplement of at most twice $0.20402$, so this
is enough to ensure that condition (2) holds as well.

Theorem~\ref{thm:extended_fuller} now tells us that
\begin{multline} \label{eqn:areas}
|\Wr(C) - \Wr(C_n)| \leq 
\frac{1}{2\pi} \sum_{i=1}^{n}
  \left| \Area R \left( T  \left( s_i \right),
                        T  \left( s_{i+1}  \right),
                        T_n\left( s_i   \right) \right) \right| \\
+ \left| \Area \triangle T   \left( s_i  \right) 
			 T_n \left( s_{i-1} \right)
                         T_n \left( s_i  \right) \right|,
\end{multline}
where the first term is the area of the spherical region bounded by
the geodesics from $T_n(s_i)$ to~$T(s_i)$ and $T(s_{i+1})$ and the
portion of $T$ between $s_i$ and $s_{i+1}$, and the second term is the
area of the spherical triangle. Our job now is to estimate the 
areas of these regions. To do so, we first invoke Taylor's Theorem,
in the form commonly used in numerical analysis (c.f. \cite{atkinson}, Thm.1.4):

\begin{theorem} {\rm (Taylor's Theorem)} Suppose $C(s)$ is a curve of class
$\mathcal{C}^4$, with fourth derivative bounded by $B'_4$. Then (choosing
coordinates so that $C(0)$ is at the origin), 
\begin{equation} \label{eqn:taylor}
C(s) = s C'(0) + \frac{s^2}{2} C''(0) + \frac{s^3}{6} C'''(0) + R_4(s),
\end{equation}
where $|R_4(s)| < s^4 B'_4$.
\end{theorem}
We will use this expression for $C(s)$ frequently in our work below.

\begin{lemma}\label{lem:length}
For any $s$, we have
\begin{equation}
|s - |C(s)|| < \frac{K^2}{24} |s^3| + \frac{1}{120} |s^5|, \qquad \text{ and } \qquad |C(s)| \leq |s|.
\end{equation}
Further, for any edge of $C_n$, the difference $|s_{i+1} - s_i|$ is at
most~$1.01\,x$.
\end{lemma}

\begin{proof}
We assume without loss of generality that $s$ is positive.  By Schur's
lemma (\cite{Chern}), since the curvature of $C$ is bounded above
by $K$, $|C(x)|$ is at least the length of a chord across an arc of
length $s$ on a circle of radius $1/K$, or~$(2/K) \sin (K/2) s$. This
means that we have
\begin{equation*}
\frac{2}{K} \sin \frac{K}{2} s = s - \frac{K^2}{24} s^3 + R_5(s),
\end{equation*}
where $R_5(s)$ is the term of order $s^5$ which comes from the usual
Taylor expansion of $\sin s$. In particular,
\begin{eqnarray*}
|s - |C(s)|| &<& s - \frac{2}{K} \sin \frac{K}{2} s \\
             &<& \frac{K^2}{24} s^3 - R_5(s),
\end{eqnarray*}
where $R_5(s) < \frac{1}{120} s^5$. The upper bound on $|C(s)|$ comes
from the fact that $C$ is unit-speed.

The second statement is another Schur's lemma calculation; this time
invoking our hypothesis that $x > (5/2) K$ and observing that $1.01
\sin y > y$ for $y$ between $0$ and $1/5$.
\end{proof}

We will also need an upper bound on the angle between $T(s)$ and $C_n(s)$.
\begin{lemma}
\label{lem:angle}
The angle between the tangent vector $T(s)$ and the corresponding
tangent vector $T_n(s)$ to $C_n$ is bounded above by
\begin{equation}
\angle T(s) T_n(s) < 0.51005 Kx.
\end{equation}
\end{lemma}

\begin{proof}
Assume that $s$ is between $s_i$ and $s_{i+1}$. Then 
\begin{equation}
\sin \angle T(s) T_n(s) = \frac{| [C(s_{i+1}) - C(s_{i})] \cross T(s) |}{|C(s_{i+1}) - C(s_{i})|}.
\end{equation}
But we have
\begin{equation*}
C(s_{i+1}) - C(s_i) = \int_{s_i}^{s_{i+1}} T(t) \, dt,
\end{equation*}
and for any $t$, we have 
\begin{equation*}
T(t) = T(s) + \int_s^t T'(u) \, du.
\end{equation*}
This means that
\begin{eqnarray}
[C(s_{i+1}) - C(s_i)] \cross T(s) &=& \int_{s_i}^{s_{i+1}} T(t) \cross T(s) \, dt \\
                                    &=& \int_{s_i}^{s_{i+1}} \int_s^t T'(u) \cross T(s) \, du \, dt.
\end{eqnarray}
Since $|T'(u) \cross T(s)| \leq |T'(u)| |T(s)| \leq \kappa(u) < K$,
and $s$ is between $s_{i}$ and $s_{i+1}$, a small computation reveals that 
this integral is bounded by $\frac{K}{2} (s_{i+1} - s_i)^2$. 

Since the length $|C(s_{i+1}) - C(s_i)|$ is bounded below by $(1/1.01)
(s_{i+1} - s_i)$ by Lemma~\ref{lem:length}, we get
\begin{equation}
\sin \angle T(t) T_n(t) < \frac{1.01}{2} K x.
\end{equation}
Since $1/x > (5/2) K$, this is always bounded above by $1.01/5$, and so 
\begin{equation}
\angle T(t) T_n(t) < \frac{1.01^2}{2} K x.
\end{equation}
\end{proof}

We are now ready to embark on the main work of the proof: estimating
the areas in Equation~\ref{eqn:areas}. We begin with the
first term: the area bounded by the portion of $T(s)$ between $s_i$
and $s_{i+1}$, together with the great circle arcs joining $T(s_i)$
and $T(s_{i+1})$ to $T_n(s_i)$. Without loss of generality, we may
assume that $i=0$, that $s_0 = 0$, and that $C(0) = \mathbf{0}$, and
apply the Taylor expansion of Equation~\ref{eqn:taylor} to $C$ at $0$.
Our strategy is to prove that this region is contained in a
neighborhood of the great circle arc joining $T(0)$ and $T(s_1)$.
Suppose $s$ is between $0$ and $s_1$. We want to bound the
height of $T(s)$ above the $T(0)$, $T(s_1)$ plane, or
\begin{equation} 
h(s) := \frac{C'(s) \cdot C'(0) \cross C'(s_1)}{|C'(0) \cross C'(s_1)|}.
\end{equation}
First, we have 
\begin{eqnarray*}
C'(s_1) &=& C'(0) + s_1 C''(0) + \frac{s_1^2}{2} C'''(0) + R_3(s_1). \\
C'(s)   &=& C'(0) + s C''(0) + \frac{s^2}{2} C'''(0) + R_3(s). 
\end{eqnarray*}
Using the triple product identities, we can rewrite $h(s)$ in terms of 
the inner product of $C'(0)$ and the cross product of these vectors.
For the triple product, we get
\begin{equation}  
\left[ \frac{s_1^2 s}{2} - \frac{s^2 s_1}{2} \right] C'(0) \cdot C''(0) \cross C'''(0) +
C'(0) \cdot \left[ R_3(s_1) \cross C'(s) + C'(s_1) \cross R_3(s) \right].
\end{equation}
Expanding the last term, we see that is the sum of a term of order
$s_1 s^3$ and a term of order $s s_1^3$. Thus, to leading order,
the norm of the entire triple product is bounded above by
\begin{equation}
|h(s)| < H := \frac{KT}{2|C'(0) \cross C'(s_1)|} s_1^3 + O(s_1^4),
\end{equation}
since $s \in [0,s_1]$. We now consider the height of $T_n(0)$ above
the $T(0)$, $T(s_1)$ plane.  Since $T_n(0)$ is the normalization of
$C(s_1) - C(0) = C(s_1)$, this height is given by
\begin{equation}
\frac{C(s_1)}{|C(s_1)|} \cdot \frac{C'(0) \cross C'(s_1)}{|C'(0)\cross C'(s_1)|}.
\end{equation}
As before, we get 
\begin{equation}
C'(0) \cross C'(s_1) = s_1 C'(0) \cross C''(0) + \frac{s_1^2}{2} C'(0) \cross C'''(0) + O(s_1^3).
\end{equation}
Taking the dot product with the Taylor expansion of $C(s_1)$, we get
only terms of order $O(s_1^4)$ and higher.  Thus, to leading order,
this region is contained in a rectangle based on the great circle arc
joining $C'(0)$ and $C'(s_1)$ of height $H$. We now estimate the
area of this rectangle.

First, we note that the length of the great circle joining $C'(0)$
and $C'(s_1)$ is given by the angle~$\theta$ between $C'(0)$ and
$C'(s_1)$. Since $s_1 < 1.01\, x$ by Lemma~\ref{lem:length}, this length
is bounded above by $1.01 \, K x$, which is less than $0.404$ by our
hypotheses on~$x$. Since $H$ is small compared to $s$, we may
assume that the entire rectangle is contained within a spherical
disk of radius $0.5$. 

We project the rectangle to the plane by central projection: this map
is increasing on lengths and areas, and increases length by at most a
factor of $1.01$. The area of the rectangle in the plane is
overestimated by the product $1.01 \, \theta H$. On the other hand, we
have $|C'(0) \cross C'(s_1)| = \sin \theta$. And for $\theta < 0.404$,
$1.02 \sin \theta > \theta$. Keeping track of the various constants
involved, and using the fact that $s_1 < 1.01 \, x$ again, the area of
this spherical region is overestimated by
\begin{equation}
\Area R(T(s_i), T(s_{i+1}), T_n(s_i)) < KT \, x^3 + O(x^4),
\end{equation}

We now turn to the second term in the Equation~\ref{eqn:areas}: the
area of the spherical triangle bounded by $T(s_i)$, $T_n(s_{i-1})$ and
$T_n(s_i)$. Without loss of generality we assume that $i=1$,
that $s_1 = 0$, and that $C(0) = \mathbf{0}$, and we expand $C$
around $0$ using Equation~\ref{eqn:taylor}.
We wish to compute
\begin{equation}
\Area \triangle\left( \frac{C(s_0)}{|C(s_0)|}, \frac{C(s_2)}{|C(s_2)|}, C'(0) \right) = 
\left| \left( \frac{C(s_0)}{|C(s_0)|} - C'(0) \right) \cross \left( \frac{C(s_2)}{|C(s_2)|} - 
C'(0) \right) \right|. 
\end{equation}
If we factor out $1/|C(s_0)||C(s_2)|$, we are left with the norm of 
the cross product of two terms:
\begin{eqnarray*}
C(s_0) - |C(s_0)| C'(0) &=& (s_0 - |C(s_0)|) \, C'(0) + \frac{s_0^2}{2} C''(0) + \frac{s_0^3}{6} C'''(0) + R_4(s_0) \\
C(s_2) - |C(s_2)| C'(0) &=& (s_2 - |C(s_2)|) \, C'(0) + \frac{s_2^2}{2} C''(0) + \frac{s_2^3}{6} C'''(0) + R_4(s_2).
\end{eqnarray*}
Using Lemma~\ref{lem:length}, we see that $|s - |C(s)|| < (K^2/24) s^3 + O(s^5)$,
and we see that the leading term of this expression contains fifth powers of 
of $s_0$ and $s_2$, and is bounded by: 
\begin{equation}
s_0^2 s_2^2 \left( \frac{K^3}{48} + \frac{KT}{12} \right) (s_0 + s_2)
\end{equation}
However, we must still divide by $|C(s_0)| |C(s_2)|$. By Lemma~\ref{lem:length}, 
we see that the ratios $s_0/|C(s_0)|$ and $s_2/|C(s_2)|$ are bounded above
by $1.01$. Thus, using the same Lemma to conclude that $s_2$ and~$s_0$ are
less than $1.01 \, x$, and making a central projection argument as before, 
we are left with
\begin{equation}
\Area \triangle (T_n(s_i), T_n(s_{i-1}), T(s_i) ) < \frac{K^3 + KT}{3} x^3 + O(x^4).
\end{equation}
Summing over $i$, and dividing by $2\pi$, then writing $K$ and $T$ in 
terms of $B_2$ and $B_3$, we obtain the statement of the theorem.
Note that we have overestimated the numerical constants to simplify
the resulting formula.
\end{proof}

If a curve has a small region of high curvature, and larger regions of
low curvature, it may be desirable to approximate the curve more
carefully in the regions of high curvature in order to save time in
the computation of writhe. Since our error bound is additive along the
curve, these methods are well suited to this case. We have

\begin{corollary} \label{cor:alternate_main}
Suppose $C$ is a $\mathcal{C}^4$ curve and $C_n$ is a curve
inscribed in $C$ so that $C$ and $C_n$ obey the hypotheses of 
Theorem~\ref{thm:disc-error}. 

If $C$ and $C_n$ are divided into regions $R_i$, each containing $n_i$
edges which are bounded above in length by $x_i$, and so that the
derivatives of $C$ are bounded by $B_{1i}$, \dots, $B_{4i}$ and 
$1/x_i > 5B_{2i}$, then
\begin{equation*}
|\Wr(C) - \Wr(C_p)| < \sum_i \alpha_i \, n_i x_i^3 + n_i O(x_i^4).
\end{equation*}
where each $\alpha_i$ is a numerical constant less than $B_{2i} (5
B_{2i}^2 + B_{3i})$.
\end{corollary}

We make one more observation:

\begin{proposition}\label{prop:non-planar}
Let $C$ be a simple, closed space curve of class $\mathcal{C}^2$, and
$C_p$ be a polygonal approximating curve as in Theorem~\ref{thm:disc-error}
or Corollary~\ref{cor:alternate_main}. 

If the arc joining the endpoints of a sequence of $n$ edges of $C_p$
is planar, then the $n-2$ edges interior to this region contribute
nothing to the error bound in the Theorem. 

In particular, this means that the derivative bounds in both statements
can be taken to be bounds on the derivatives of the {\em non-planar}
regions of the curve $C$.
\end{proposition}

\begin{proof}
On these edges, the tantrix of the smooth curve and the polygonal curve
parametrize the same great circle arc on $S^2$. Thus, the ribbon between
these curves has zero area.
\end{proof}
 
\section{Example Computations \label{sec:examples}}

We are now prepared to test Theorem~\ref{thm:main} by computing the
writhing numbers of various polygonal approximations of a smooth
curve, and comparing the results to the exact writhe of the smooth
curve.  To control the numerical error introduced in these
calculations, all of these computations were performed using an
arbitrary-precision implementation of Banchoff's formula for the
writhing number of a polygonal curve. The initial runs were performed
with $45$ decimal digits of precision. They were checked against runs
performed with $54$ digits of precision. Since the results agreed,
we feel confident that roundoff error does not affect the computations
reported on below.

The curve whose writhe we computed is an example of Fuller\cite{fuller}:
\begin{center}
\begin{tabular}{cc}
\includegraphics[height=1.5in]{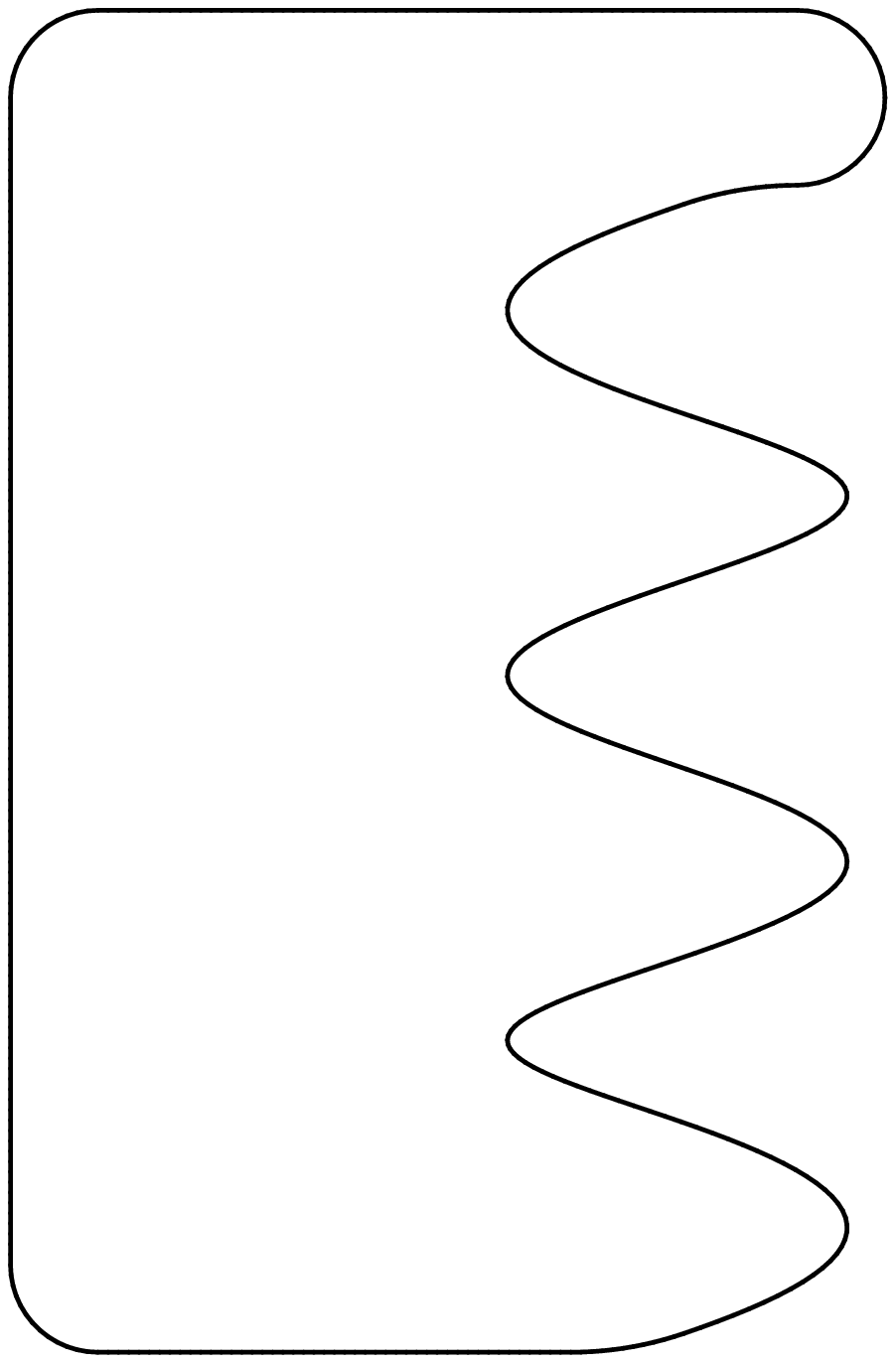} & 
\includegraphics[height=1.5in]{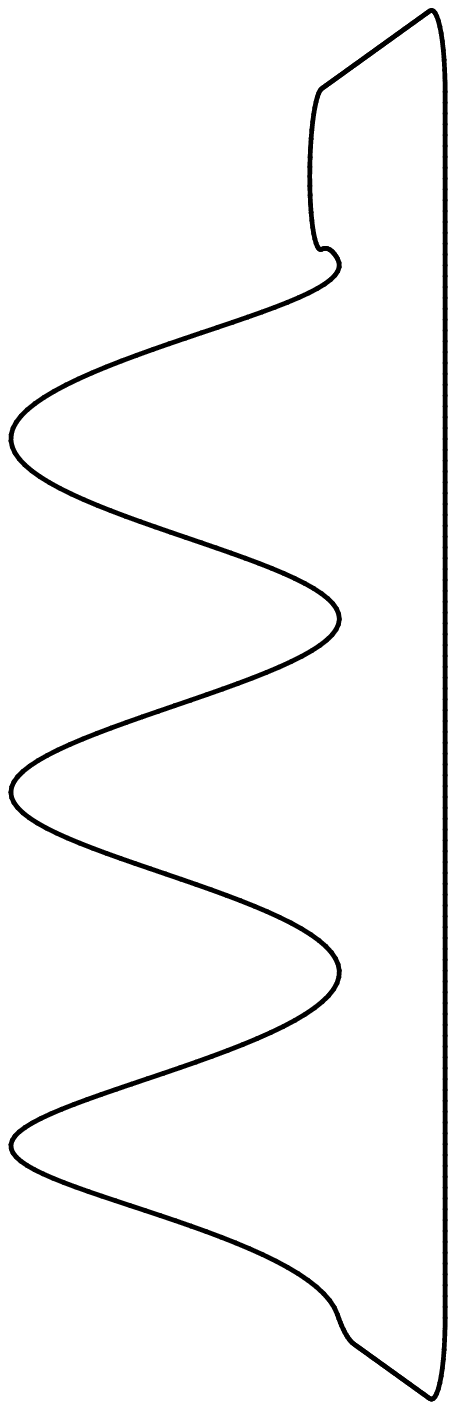} \\
\end{tabular}
\end{center}

\begin{center}
\parbox{5in}{
\noindent{\bf Figure 4.} This example of Fuller's ``closed helix'' is
composed of $3$ turns of a helix of radius $1$ with pitch angle
$0.33$, with ends joined by a planar curve. }
\end{center}
\medskip

Using Theorem~\ref{thm:spherical_area}, and the \Calugareanu-White
formula, it is easy to see that the writhe of this curve is $3(1 -
\sin 0.33) \simeq 2.0278709$.  After all, the area enclosed by the
tantrix of this curve $C$ is that of a hemisphere, plus $3$ enclosures
of a spherical cap of radius~$\pi/2 - 0.33$. Thus the writhe of the
curve is equal to $1-\sin 0.33 \mod 2$. To complete the computation,
one sets up a frame on the curve, and computes its twist and linking
number.  (Details for this computation can be found
in~\cite{fuller}.)

We now take a series of polygonal approximations to $C$, and compare
the difference between their writhing numbers and the writhe of $C$
to the bounds of Theorem~\ref{thm:disc-error}.

We begin by finding bounds on the derivatives of $C$ and the edge
length of our approximations. By Proposition~\ref{prop:non-planar}, it
suffices to find derivative bounds for the helical region of $C$.
Since the helix has unit radius, both $B_2$ and $B_3$ can be taken to
be one. The curve is parametrized so that $|C'(s)| \geq 1$.

Here are the results of computing writhe with various numbers of edges: \\

\begin{center}
\begin{tabular}{r|c|c|c|c}
$n$ & $\Wr(C_n)$ & $|\Wr(C_n) - \Wr(C)|$ & $x$ & $\alpha \, nx^3$ \\
\hline
100   & $2.00541$ & $0.02246$ & $0.506$   &  77.73  \\  
250   & $2.02434$ & $0.00353$ & $0.203$   &  12.55  \\ 
500   & $2.02697$ & $0.0009$  & $0.101$   &  3.09  \\ 
1000  & $2.02763$ & $0.00024$ & $0.051$   &  0.786  \\
\end{tabular} 
\end{center}
\medskip

It is worth examing a graph of these results.
\smallskip
\begin{center}
\makebox[5in][l]{
\includegraphics{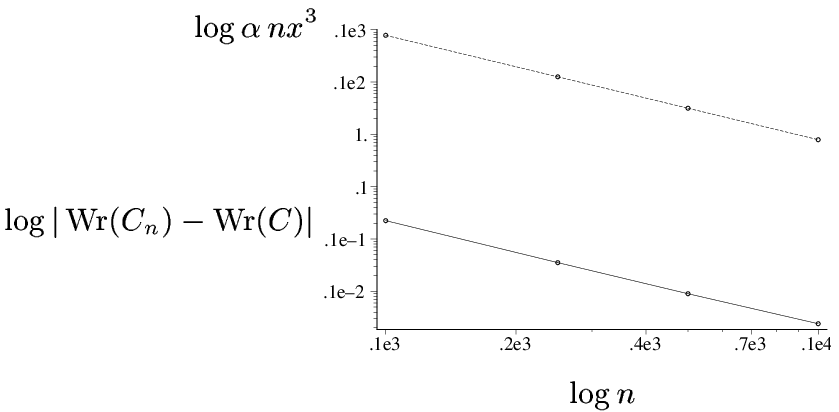} \hspace{1in}
} \\
\end{center}
\smallskip
\begin{center}
\parbox{5in}{ {\bf Figure 5.} This graph shows a log-log plot of the
actual error in computing the writhing number for one of Fuller's
``closed helices'' with various numbers of edges (lower solid line),
together with our error bounds (upper dotted line). The fact that the lines
are parallel shows that the convergence is of order $n^2$, as
predicted by Theorem~\ref{thm:disc-error}. }
\end{center}
\medskip

\section{Further Directions}

In this paper, we have given a set of asymptotic error bounds which
allow us to compute the writhe of a closed space curve with defined
accuracy by computing the writhe of a polygonal approximation to this
curve. The example we computed in Section~\ref{sec:examples} shows
that our bounds are of the right order of magnitude: roughly speaking,
the writhe converges quadratically in the number of edges of the
approximation. Our work leaves open several directions for further
inquiry. 

First, it is puzzling that our approximation theorem for curves with
corners (Proposition~\ref{prop:approx}) should depend on the
hypothesis that the limit curve is polygonal. To sketch an extension
of this result, we recall a definition from Chern (\cite{Chern}):
\begin{definition} 
\label{def:tantrix}
The tantrix of a piecewise $C^1$ curve $C(s)$ with positive
corner angles is the image of $T(s)$ on the unit sphere, together with
the great circle arcs joining the pairs of tangent vectors at each
corner of the curve.
\end{definition}

We note that our Theorem~\ref{thm:extended_fuller} shows that Fuller's
$\Delta\!\Wr$ formula holds for polygonal curves with the
definition of tantrix extended as above. We further suspect that:

\begin{conjecture}
Fuller's Spherical Area formula (Theorem~\ref{thm:spherical_area}) and
Fuller's $\Delta\!\Wr$ formula (Theorem~\ref{thm:fuller_delta_writhe})
hold for piecewise $C^2$ curves with the extended definition of
tantrix given by Definition~\ref{def:tantrix}.
\end{conjecture}

The proofs of both of these theorems depend on the \Calugareanu-White
formula, which only applies to closed curves. Thus all of our results
are restricted to closed curves. This leaves open a much more
important problem:
\begin{problem}
Extend all these theorems (the \Calugareanu-White
formula, Fuller's spherical area formula, and Fuller's $\Delta\!\Wr$
formula) to open curves.
\end{problem}
In particular, extending the results of this paper to open curves
would be useful for applications in biology, where the curves of interest
are not neccesarily closed. We note that while Fuller's
$\Delta\!\Wr$ formula makes sense for open curves, computational
examples show that it does not give the correct answer: boundary terms
must be added to account for the ends of the curves.

\section{Acknowledgements}

I am grateful to Herbert Edelsbrunner, Herman Gluck, Issac Klapper,
John Maddocks, Kathleen Rogers, and David Swigon, among many others,
for fruitful conversations. This work was supported by both a National
Defense Science and Engineering Graduate Fellowship and an NSF
Postdoctoral Research Fellowship (DMS-99-02397).

\bibliography{thick}
\bibliographystyle{plain}
 
\end{document}